\documentclass[12pt,amssymb,amsmath]{article}
\usepackage{epsf}
\usepackage{amssymb}
\newcommand{\text}{\rm}

\newcommand{\beqn}{\begin{eqnarray}}
\newcommand{\eeqn}{\end{eqnarray}}

\newcommand{\beq}{\begin{equation}}
\newcommand{\eeq}{\end{equation}}
\newtheorem{defn}[subsection]{Definition}
\newtheorem{thm}[subsection]{Theorem}
\newtheorem{prop}[subsection]{Proposition}
\newtheorem{cor}[subsection]{Corollary}

\newenvironment{rem}{\smallskip\noindent%
\refstepcounter{subsection}%
{\bf \thesubsection}~~{\sc Remark.}\hspace{-1mm}}
{\smallskip}

\newenvironment{ex}{\smallskip\noindent%
\refstepcounter{subsection}%
{\bf \thesubsection}~~{\sc Example.}}
{\smallskip}

\newenvironment{que}{\smallskip\noindent%
\refstepcounter{subsection}%
{\bf \thesubsection}~~{\sc Question.}\hspace{-1mm}}
{\smallskip}

\begin{document}

%\hskip 2.3in{}

\title{\Large\bf  Euler equations on homogeneous spaces and
Virasoro orbits}
\date{December 14, 2001 \\~\\To appear in {\it Adv. Math.}}
\author{\large\bf Boris Khesin\thanks{Department of Mathematics,
University of Toronto, Toronto, ON M5S 3G3, Canada;
e-mail: {\tt khesin@math.toronto.edu}}~
and Gerard Misio\l ek\thanks{Department of Mathematics,
University of Notre Dame, Notre Dame, IN 46556, USA;
e-mail: {\tt misiolek.1@nd.edu}}~}
\maketitle
\bigskip

\begin{abstract}
We show that the following three systems related to various 
hydrodynamical approximations: the Korteweg--de Vries equation,
the Camassa--Holm equation, and the Hunter--Saxton equation,
have the same symmetry group and similar bihamiltonian structures. 
It turns out that their configuration space is the Virasoro 
group and all three dynamical systems  can be regarded as equations of the
geodesic flow associated to different right-invariant metrics on this 
group or on  appropriate homogeneous spaces. 
In particular, we describe how  Arnold's approach to the Euler 
equations as geodesic flows of one-sided invariant metrics extends from 
Lie groups to homogeneous spaces.

We also show that the above three cases describe 
all generic bihamiltonian systems which are
related to the Virasoro group and can be integrated by the translation
argument principle: they correspond 
precisely to the
three different types
of generic Virasoro orbits. Finally, we discuss 
interrelation between the above metrics  and
Kahler structures on Virasoro orbits as well as open questions 
regarding integrable systems corresponding to a
finer classification of the orbits.
\end{abstract}

%%%%%%%%%%%%%%%%%%%%%%%%%%%%%%%%%%%%%%%%%%%%%%%%%%%%%%%%%%%%%

\section{Introduction}

\setcounter{equation}{0}

One of the main mechanisms of integrability of evolution  equations
is the presence of two compatible Hamiltonian structures.
In this paper we compare Hamiltonian  properties of  three extensively
studied nonlinear equations of mathematical physics, related to
various hydrodynamical approximations: the Korteweg-de Vries equation
\begin{equation}
u_t=-3uu_x+cu_{xxx},\label{eq:KDV}
\end{equation}
the Camassa-Holm equation
\begin{equation}
 u_t -  u_{txx} =
- 3u u_x+ 2 u_x  u_{xx} +u u_{xxx} +c  u_{xxx}
\label{eq:CH}
\end{equation}
derived as a shallow water equation in \cite{ch} (see also the paper
\cite{ff}),
and the Hunter--Saxon equation \cite{hs}
\begin{equation}
u_{txx}=-2 u_x  u_{xx} -u u_{xxx},\label{eq:HS}
\end{equation}
describing weakly nonlinear unidirectional waves.
All three equations are known to be bihamiltonian and to
possess infinitely many conserved quantities,
as well as remarkable soliton or soliton--like solutions.
A motivation for our paper was the paper \cite{bss}, which described
scattering theory for all three equations in a unified way.

As  we show in this paper,  the main reason why such a common treatment is
possible is that all these equations have the same symmetry group.
It turns out to be the Virasoro group,
a one-dimensional extension of the group
of smooth transformations of the circle.
More precisely, the Virasoro group serves as the configuration space, and
all three equations can be regarded as  equations of the geodesic flow
related to different right-invariant metrics on this group (in the case of
KdV and CH) or on an associated homogeneous space (in the  case of HS).
(Here we will be mostly concerned with the periodic case,
though many statements can be extended to the case of rapidly
decaying potentials on the real line.)

The main goal of this paper is to give a  description of the three
equations as {\it bihamiltonian} systems on the  dual
to the Virasoro algebra, and to
relate them to the geometry of the Virasoro coadjoint orbits.
One of the corresponding Hamiltonian (or Poisson) structures is provided
by the linear Lie-Poisson bracket, and it is the same for all three equations.
The other Poisson structure is constant and can be viewed as
a ``linear structure frozen at a point.'' The corresponding
``freezing points'' are different for each  equation. We will also
see that, in a sense,  Equations (\ref{eq:KDV}), (\ref{eq:CH}), and
(\ref{eq:HS}) exhaust all generic possibilities, and among them
the Camassa-Holm equation (\ref{eq:CH})
is the ``most general'' equation that can be
obtained  by the ``freezing argument''  method on the dual Virasoro space.

We tried to make the paper self-contained, including in
it necessary background
on the Euler equations and the classification  of Virasoro orbits.
For additional information
we refer the interested reader to the expositions in
\cite{seg} or \cite{ak}, as well as to the  original papers listed in
the bibliography.

\medskip

In \cite{arn} V. Arnold suggested a general framework for the Euler
equations on an arbitrary (possibly infinite-dimensional) group, which we
recall below. In this framework the Euler
equation describes a geodesic flow with respect to a suitable
one-sided invariant Riemannian metric on the given group.

In  Section 2 we show how Arnold's approach to the Euler
equation works for the Virasoro group and provides a natural geometric setting
for the Korteweg--de Vries and Camassa--Holm equations. In Section 3
we give a Hamiltonian reformulation of the Euler equation.
In  Section 4 we extend this approach to include
geodesic flows on homogeneous spaces and then use it to describe
the Hunter--Saxton equation and its relatives.\footnote{In particular,
the Harry Dym  equation \cite{kru}, \cite{hzh} can be found as one
of the equations in the bihamiltonian hierarchy associated with the HS
system. This equation was also  considered
in \cite{bss}. In this way, the Harry Dym equation also becomes associated
to the geodesic interpretation.}
This extension might be thought of as a version of the Hamiltonian formalism 
for homogeneous spaces developed in \cite{th}, \cite{gs}, which is applied 
to the  case of a degenerate metric and in our infinite-dimensional situation.

In Sections 5--6 we develop the  bihamiltonian formalism for the
Euler systems. We show why the  three equations above represent
three main classes (and exhaust the rotation-invariant)
bihamiltonian systems on the Virasoro algebra
that can be integrated by means of the ``freezing'' (called also, the 
translation of) argument method.
It turns out that the above three  equations correspond to three
 different types of Virasoro coadjoint orbits of low codimensions.

In our classification of Poisson pairs and the corresponding equations
we rely heavily on the  classical classification
of the Virasoro  orbits, and we recall it in the Appendix.
An interesting open  question is to extend the classification
of the equations to orbits of higher codimension, as well as to
show how the discrete invariant of Virasoro orbits manifests itself
in the related bihamiltonian systems.

%%%%%%%%%%%%%%%%%%%%%%%%%%%%%%%%%%%%%%%%%%%%%%%%%%%%%%%%%%%%%%
\bigskip

%%%%%%%%%%%%%%%%%%%%%%%%%%%%%%%%%%%%%%%%%%%%%%%%%%%%%%%%%%%%%%%%%%%%%
\section{The Euler properties of the KdV, CH, and HS equations}

The main objects in our consideration  will be the
group $\mathrm{Diff}(S^1)$ of all diffeomorphisms
of a circle, its Lie algebra $vect(S^1)$ of vector fields,
and their central extensions.
The following nontrivial one-dimensional
extension of the algebra of vector fields has a special name.

\begin{defn}
The {\rm Virasoro algebra} is an extension of the
Lie algebra ${\text{vect}}(S^1)$ of vector fields on the circle:
$$
{vir}
=
vect(S^1) \oplus {\mathbb{R}}
$$
with the commutator between pairs (consisting of a vector field and a real
number) given by
$$
[(v(x)\partial_x,b),(w(x)\partial_x,c)]
=
 \left(~ (-vw_x + v_x w)(x)\partial_x , \, C(v\partial_x,w\partial_x)~ \right)
$$
where
$$
C(v\partial_x,w\partial_x)
=
\int_{S^1}  vw_{xxx} \, dx
$$
is the {\rm Gelfand--Fuchs cocycle.} (Here $x$ is a coordinate on
the circle, the subscript $x$  stands for the  derivative in $x$, and
$\partial_x$ denotes the vector field
$\frac{\partial}{\partial x}$ on $S^1$.)
\end{defn}

Define the following two-parameter family of quadratic forms, 
``$H^1_{\alpha, \beta}$-energies,'' on the Lie algebra ${vir}$: 
%by means of the $L^2$ and $H^1$ inner products:
\begin{equation}
\langle (v(x)\partial_x,b) , (w(x)\partial_x,c) \rangle_{H^1_{\alpha, \beta}}=
\int_{S^1} (\alpha\, vw + \beta\, v_x w_x) \, dx + bc\,.
\label{eq:metric}
\end{equation}
The case $\alpha=1, ~~ \beta=0$ corresponds to the  $L^2$ inner product,
while  $\alpha=\beta=1$ corresponds to $H^1$. 
Given  $\alpha\not= 0$ and any $\beta$, extend the $H^1_{\alpha, \beta}$-energy
 to a right-invariant metric
on the Virasoro group $Vir$. This group corresponds to the Virasoro algebra
and is defined as follows.

\begin{defn}
The {\rm Virasoro group} $Vir$
is  the   product
$$
Vir=\mathrm{Diff}(S^1) \times {\mathbb{R}},
$$
where the group multiplication between the pairs is given by
$$
(\psi(x) , a) \circ (\phi(x) , b)
=
(~ (\psi \circ \phi)(x) , \, a + b + B(\psi,\phi)~)
$$
and
$$
B(\psi,\phi):=\int_{S^1}
\log ((\psi \circ \phi)_x) \, d\log (\phi_x)
$$
is  the {\rm Bott cocycle.}
(Here the diffeomorphisms
of $S^1$ are described by functions, e.g., $x\mapsto \phi(x)$.)
\end{defn}

Having equipped the Virasoro group with those right-invariant metrics,
one can consider the geodesic flows they generate.

\begin{thm}\label{Euler} 
1) {\rm{\cite{ok}}}
The KdV equation is the Euler equation, describing
the geodesic flow on the Virasoro group with respect to
the right-invariant $L^2$-metric.

2) {\rm{\cite{mi}}} 
The CH equation is the Euler equation for
the geodesic flow on the same group with respect to
the right-invariant Sobolev $H^1$-metric.
 
\end{thm}

It turns out that one can give a similar description of the Hunter-Saxton
equation as a geodesic flow on a {\it homogeneous space} related to
the Virasoro algebra. Consider the 
$\dot{H}^1$ -quadratic 
form (which is the  $H^1_{\alpha, \beta}$ -form with  
$\alpha=0, ~~ \beta=1$) on the Virasoro algebra:
\begin{equation}
\langle (v(x)\partial_x,b) , (w(x)\partial_x,c) 
\rangle_{\dot{H}^1}=
\int_{S^1}  v_x w_x\, dx + bc\,.
\label{eq:metric0}
\end{equation}
Although this form is degenerate, as is the corresponding right-invariant 
metric on the Virasoro group, one can  define a nondegenerate metric by
descending on an appropriate quotient space.

\begin{thm}\label{homEuler} The HS equation is the equation describing
the geodesic flow on the homogeneous space
$Vir/\mathrm{Rot}(S^1)$ of the Virasoro group modulo rotations
with respect to the right-invariant homogeneous
$\dot{H}^1$ metric.
\end{thm}

Note that one can also obtain the HS equation by considering
the smaller homogeneous space
$\mathrm{Diff}(S^1)/\mathrm{Rot}(S^1)$ of all diffeomorphisms of the circle
modulo rotations, as we explain below.

\smallskip

These three equations  essentially exhaust the list of
integrable systems associated with the Virasoro algebra
and integrated by the freezing argument method, as we discuss
below. Note that their degenerations include, e.g., the
inviscid Burgers equation (it corresponds to the $L^2$-metric on 
the ``centerless Virasoro''  group,
$\mathrm{Diff}(S^1)$).

\begin{rem} Before proving the theorems,  we  recall
the general set-up for the Euler equation on an arbitrary
Lie group, suggested by  V.~Arnold in \cite{arn}.
Consider a (possibly infinite-dimensional) Lie group $G$,
which can be thought of as the configuration
space of some physical system. (Examples from \cite{arn}:
$SO(3)$ for a rigid body
or  the group $\mathrm{SDiff}(M)$ of volume-preserving diffeomorphisms
for an ideal fluid filling a domain $M$.)
The tangent space at the identity of the
Lie group $G$ is the corresponding Lie algebra
$\mathfrak g$.
Fix some (positive definite) quadratic form, the ``energy,''
on  $\mathfrak g$. We consider right translations
of this quadratic form to the tangent space at any point of the group
(the ``translational symmetry'' of the energy).
This way the energy defines a  right-invariant
Riemannian metric on the group $G$.
The geodesic flow on $G$ with respect to this energy metric
represents the extremals of the least action principle, i.e.
the actual motions of our physical system.\footnote{For a rigid body one
has to consider  left translations, but in our exposition we stick
to the right-invariant case in view of its applications to the groups
of diffeomorphisms.}

To describe a geodesic on the {\it Lie group} with an initial velocity $v(0)$,
we transport its velocity vector at any moment  $t$
to the identity of the group (by using the right translation).
This way we obtain the evolution law for $v(t)$, given by
a (non-linear) dynamical system $dv/dt=F(v)$
on the {\it Lie algebra} $\mathfrak g$.

\begin{defn} The system on the Lie algebra $\mathfrak g$,
describing the evolution
of the velocity vector  along a geodesic in a right-invariant metric
on the Lie group $G$, is called the {\rm Euler equation} corresponding to
this metric on $G$.
\end{defn}

%%%%%%%%%%%%%%%%%%%%%%%%%%%%%%%%%%%%%%%%%%%%%%%%
%%%%%%%%%%%%%%%%%%%%%%%%%%%%%%%%%%%%%%%%%%%%%%
\begin{figure}[htb]
\centerline{\epsfbox{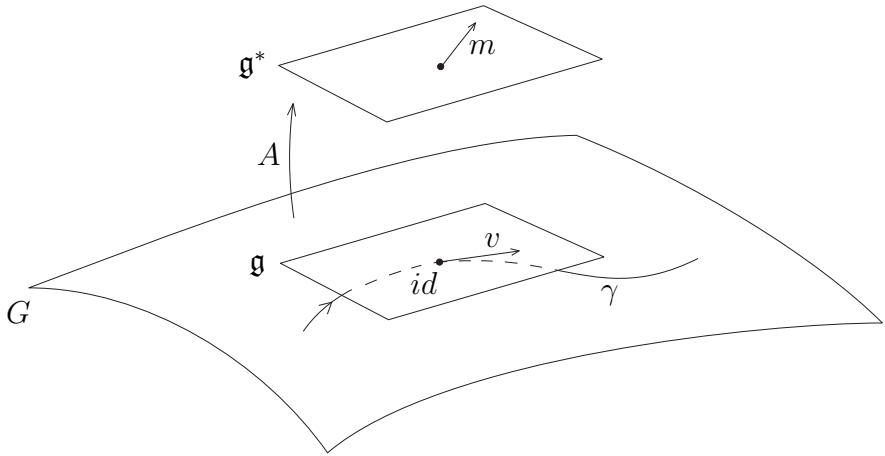}}
\begin{picture} (0, 64)
\put(143,210){${\mathfrak g}^*$}
\put(230,217){$m$}
\put(150,175){$A$}
\put(147,135){${\mathfrak g}$}
\put(236,144){$v$}
\put(208,126){$id$}
\put(280,124){$\gamma$}
\put(55,115){$G$}
\end{picture}
\vskip -24mm \caption{\small \footnotesize
The vector $v$ in the Lie algebra $\mathfrak g$
traces the evolution of the velocity vector of a geodesic $\gamma$
on the group. The inertia operator $A$ sends $v$ to  a vector $m$
in the dual space $\mathfrak g^*$.} \label{fig1}
\end{figure}
%%%%%%%%%%%%%%%%%%%%%%%%%%%%%%%%%%%%%%%%%%%%
%%%%%%%%%%%%%%%%%%%%%%%%%%%%%%%%%%%%%%%%%%%%%%%%%%

In particular,  the above scheme works for the Virasoro group (see Theorem
\ref{Euler}) and allows one to describe the Korteweg--de Vries
and Camassa--Holm equations as geodesic equations on that group.
It also  can be extended to include
geodesic flows on homogeneous spaces and to describe
the Hunter--Saxton equation, as we discuss below.
\end{rem}

\bigskip

%%%%%%%%%%%%%%%%%%%%%%%%%%%%%%%%%%%%%%%%%%%%%%%%%%%%%%%%

\section{Hamiltonian framework for the Euler equations}

We start with  preliminaries on Lie algebras and
Poisson structures.

\begin{defn}
The dual space $\mathfrak g^*$ to any Lie algebra $\mathfrak g$
carries a natural {\rm Lie--Poisson structure}:
$$
\{f,g\}_{LP}(m):=\langle [df,dg], m\rangle
$$
for any $m\in \mathfrak g^*$ and any two smooth functions 
$f,g$ on $\mathfrak g^*$.
(Here the differentials are taken at the point $m$, and
$\langle \cdot,\cdot\rangle$ is a natural pairing between $\mathfrak g$
and $\mathfrak g^*$.)
\end{defn}

In other words, the Lie--Poisson bracket of two linear functions on 
$\mathfrak g^*$
is equal to their commutator as elements of the Lie algebra 
$\mathfrak g$ itself.

\begin{prop}\label{hamfield}
The {\it Hamiltonian vector field} on $\mathfrak g^*$ corresponding
to a Hamiltonian function $f$ and computed  with respect to the
Lie-Poisson structure has the following form:
\begin{equation}\label{eq:hamLP}
\frac{dm}{dt}={\rm ad}^*_{df}m.
\end{equation}
\end{prop}

{\bf Proof.}
Let $dm/dt=X_f$ be the corresponding Hamiltonian  field. Then
for any function $g\in C^\infty(\mathfrak g^*)$ one has the identities
$$
i_{X_f} dg|_m= L_{X_f} g|_m=\{f,g\}_{LP}(m)=\langle [df,dg], m\rangle=
\langle dg, {\rm ad}^*_{df}m\rangle.
$$
This implies that $X_f={\rm ad}^*_{df}m$.~~$\square$

\smallskip

\begin{rem}
The differential-geometric description of the Euler equation
as a geodesic flow on a Lie group has a Hamiltonian reformulation.

Fix the notation $E(v)=\frac 12\langle v,Av\rangle$
for the energy quadratic form on $\mathfrak g$ which we used to define the
Riemannian metric. Identify the Lie algebra and its dual
with the help of this quadratic form.
This identification   $A:\mathfrak g\to \mathfrak g^*$ (called the {\it inertia
operator}) allows one to rewrite the Euler equation on the dual space
 $\mathfrak g^*$, see Fig.1.

It turns out that the Euler equation on $\mathfrak g^*$
is Hamiltonian with respect to the  Lie--Poisson structure \cite{arn}.
Moreover, the corresponding Hamiltonian function is minus the energy
quadratic form lifted from the Lie algebra to its
dual space by the same identification:
$-H(m)=-\frac 12\langle A^{-1}m,m\rangle$, where $m=Av$.
Here we are going to take it as the {\it definition}
of the Euler equation (we use the Proposition above and the  observation
$dH(m)=A^{-1}m$).
\end{rem}

\begin{defn}\label{linEuler}
The  {\rm Euler equation on} $\mathfrak g^*$ corresponding to the Hamiltonian
$-H(m)=-\frac 12\langle A^{-1}m,m\rangle$ is given by the following
explicit formula:
$$
\frac{dm}{dt}=-{\rm ad}^*_{A^{-1}m}m,
$$
as an evolution of a point $m\in \mathfrak g^*$.
\end{defn}

\begin{rem}
The underlying reason for
the Hamiltonian reformulation is the fact that any
geodesic problem in {\it Riemannian} geometry
can be described in terms of {\it symplectic} geometry.
Geodesics on $M$ are extremals of a quadratic Lagrangian (metric)
on $TM$. They  can also be described by the Hamiltonian flow on $T^*M$
for  the quadratic Hamiltonian function obtained from 
the Lagrangian via the Legendre transform.

If the manifold is a group $G$ with a right-invariant
metric then there exists the group action on the tangent bundle $TG$, 
as well as on the cotangent bundle $T^*G$.
By taking the quotient with respect to the group action, we obtain
from the (symplectic) cotangent bundle $T^*G$  the Lie-Poisson structure
on the cotangent space $T^*G|_e=\mathfrak g^*$, i.e., on the dual to
the Lie algebra.
The Hamiltonian function on $T^*G$ is dual to the Riemannian
metric (viewed as a form on $TG$),
and its restriction to $\mathfrak g^*$ is
 the quadratic form  $H(m)=\frac 12\langle A^{-1}m,m\rangle$,
~ $m\in \mathfrak g^*$.

The geodesics of a left-invariant metric on $G$ correspond to the
Hamiltonian function $H(m)$, while those of a right-invariant metric
correspond to $-H(m)$.
\end{rem}

\bigskip
Now we are ready to prove Theorem \ref{Euler}
on the Eulerian nature of the KdV and CH equations in the following 
slightly more general setting.

\begin{thm}{{\bf (=\ref{Euler}$'$})}\label{Euler'} 
The Euler equation describing
the geodesic flow on the Virasoro group with respect to
the right-invariant $H^1_{\alpha, \beta}$ -metric with $\alpha\not=0$ 
has the form:
\begin{equation}\label{abEuler}
\alpha(v_t+3vv_x)-\beta(v_{xxt}+2v_xv_{xx}+vv_{xxx})-bv_{xxx}=0,
\end{equation}
$$
b_t=0.
$$
\end{thm}

\begin{rem} By choosing  $\alpha=1, ~\beta=0$ one obtains the KdV equation,
related to the $L^2$-metric on the Virasoro algebra \cite{ok}. Similarly,
for $\alpha=\beta=1$ one recovers
a general form of the CH equation \cite{mi}. Note that by shifting
$v\mapsto v+$const we get another form of the CH equation,
in which the term $v_{xxx}$ is replaced by  $v_x$.
Finally, if $\alpha=0,~\beta=1$ then Equation (\ref{abEuler}) becomes the
HS (Hunter-Saxton) equation, which we discuss in the next section.

The case $b=0$ corresponds to considering
the non-extended Lie algebra $vect(S^1)$ of vector fields rather than the
Virasoro algebra $vir$. Depending on the values of $\alpha $ and $\beta$
one obtains either the inviscid Burgers (also called, Hopf) equation
$ v_t+3vv_x=0$ or the non-extended CH equation.
\end{rem}

{\bf{Proof of Theorem \ref{Euler'}.}} Recall that the Virasoro
coadjoint action can be computed as follows. Let
$\{(u(dx)^2,\,a)|~u\in C^\infty(S),\, a\in \mathbb R\}$ be the dual 
space to the Virasoro algebra with the natural pairing given by
$$
\langle (u(dx)^2,\,a) , (w\partial_x,c) \rangle=\int_{S^1} uw~dx+ac.
$$
(Here we denote by $u(x)(dx)^2$ (or, shorter, by $u(dx)^2$) a quadratic 
differential on the circle.)
The coadjoint operator is defined by the identity
$$
\langle {\rm ad}^*_{(v\partial_x,b)}(u(dx)^2,\,a) , (w\partial_x,c) \rangle
=
\langle (u(dx)^2,\,a) , [(v\partial_x,b) , (w\partial_x,c)] \rangle.
$$
Using the definition of the Virasoro commutator and  integrating by parts
we obtain that the right-hand-side is equal to
$$
\int_{S^1} w(2uv_x+u_xv-av_{xxx})~dx.
$$
Thus the coadjoint operator is
\begin{equation}\label{virad*}
{\rm ad}^*_{(v\partial_x,b)}(u(dx)^2,\,a)=((2uv_x+u_xv-av_{xxx})(dx)^2,\,~0).
\end{equation}

\smallskip

The $H^1_{\alpha, \beta}$ -energy (\ref{eq:metric}) on the Virasoro algebra
$$
\langle (v\partial_x,b), (w\partial_x,c) \rangle
=\int_{S^1}(\alpha \,vw+\beta \,v_x w_x)~dx+bc=\int_{S^1} v\,\Lambda w~dx+bc.
$$
corresponds to the  general inertia operator $A:vir\to vir^*$, given by
$$
(v\partial_x,b)\mapsto ((\Lambda v)(dx)^2,\, b),
$$
where $\Lambda:=\alpha -\beta \partial_x^2$ is a second order
differential operator.
This operator is non-degenerate on $vir$ for $\alpha\not=0$, while for
$\alpha =0$ it  has a non-trivial
kernel consisting of constant vector fields on $S^1$.

Now the Euler equation 
$$
\frac{d}{dt}(u(dx)^2,a)=-{\rm ad}^*_{A^{-1}(u(dx)^2,a)}(u(dx)^2,a)
$$
on $vir^*$ (see Definition (\ref{linEuler})) assumes the form
$$
\frac{d}{dt}(u(dx)^2,a)=
-\left((2u\Lambda^{-1}u_x+u_x\Lambda^{-1}u-
a\Lambda^{-1}u_{xxx})(dx)^2,~0\right).
$$
(Here we substituted $(v\partial_x, b)=A^{-1}(u(dx)^2,a)=
((\Lambda^{-1}u)(dx)^2,a)$  into the expression for $ad^*$.)

In terms of $v=\Lambda^{-1}u$  (the first component of) this equation becomes
$$
\frac{d}{dt}(\Lambda v)=-2(\Lambda v)v_x-(\Lambda v_x)v+bv_{xxx},
$$
which is equivalent to the equation (\ref{abEuler}),
since $\Lambda=\alpha -\beta \partial_x^2$.
For the second component we find that $b$ does not change in time:
$b_t=0$.~~$\square$

\begin{rem}\label{image}
In the proof we  assumed that the inertia operator
is invertible. In the next section we discuss the
precise relation of the geodesic and Hamiltonian approaches in
the  case of a degenerate metric and 
show what reductions are
necessary  for the corresponding Euler equation to make sense.
\end{rem}

\bigskip

%%%%%%%%%%%%%%%%%%%%%%%%%%%%%%%%%%%%%%%%%%%%%%%%%%%%%%%%

\section{The Euler equations on homogeneous spaces}\label{homsect}

Let $G$ be a Lie group and $K$ its subgroup. Consider the space
$G/K$ of {\it right} cosets $\{Kg~|~g\in G\}$.
Then the group $G$ acts on them on the  right.
Here we are going to develop the formalism for the Euler equation,
describing the geodesic flow on $G/K$ with respect
to a right-invariant metric.

One immediately encounters the following difficulty:
not every {\it right-invariant} metric on the group $G$, degenerate along $K$
at the identity,
descends to a metric on the space of right cosets $G/K$
(see Example (\ref{exs}a) below).
To formulate the condition, which the degenerate metric should satisfy,
 let us consider the corresponding problem  at the level of Lie algebras.

Let $\mathfrak g$ be a Lie algebra, and $A: \mathfrak g\to\mathfrak g^*$ a degenerate
inertia operator. Suppose that the kernel of $A$ is a
Lie subalgebra $\mathfrak k$. (In other words, the corresponding energy form
$E(v)=\frac 12\langle v, Av\rangle $ vanishes for all $v\in \mathfrak k$.)
Consider the right-invariant degenerate metric $E_G$ on the group $G$
obtained by translating the quadratic form
$E$ from identity to any point of the group.

\begin{thm}\label{coset} The right-invariant form $E_G$ on a group  $G$
descends to a form on the  space $G/K$ of right cosets
if and only if the quadratic form $E$ on the Lie algebra $\mathfrak g$
vanishes on the Lie subalgebra $\mathfrak k$ and  is
$\rm{Ad}$-invariant with respect to the action
of this subalgebra.
\end{thm}

\begin{rem} The condition of \rm{Ad}-invariance for $E$ reads as follows:
$$
\langle \mathrm{ad}_w v, Au\rangle =-\langle v, A(\mathrm{ad}_w u)\rangle
$$
for all $w\in \mathfrak k$ and any  $u,v\in \mathfrak g$.
The above is an infinitesimal version of the invariance of $E$ with respect to
the subgroup action:
$$
\langle (\mathrm{Ad}_k v), A(\mathrm{Ad}_k u)\rangle=\langle  v, A u\rangle
$$ for all $k\in K$.
\end{rem}

\begin{ex}\label{exs} a) (Rotations of a rod.)
The configuration space of a rod  in  ${\mathbb R}^3$
fixed at its  center of mass is  $S^2$.
It can be obtained from the configuration space of a rigid body
by moding out rotations about one of its axes: $S^2=S^1\setminus SO(3)$.
\end{ex}

Suppose that $A: {\mathbb R}^3\to {\mathbb R}^3$ is a degenerate inertia operator
 with one vanishing eigenvalue.
The corresponding eigenvector generates the 1-dimensional rotation subgroup
$S^1$ in $SO(3)$. It is not difficult to see
that the bi-invariance condition imposes the  following restriction on
$A$: its non-vanishing  eigenvalues must be equal. (Indeed,
the $S^1$-action sends one of these two eigenvectors to the other.)
Then the corresponding degenerate metric on  $SO(3)$
descends to  (a multiple of) the standard metric on the sphere,
which is the space of cosets:
$S^2=S^1\setminus SO(3)$. Geodesics with respect to the standard metric on
$S^2$ are the great circles. These geodesics describe all free
motions of the rod.\footnote{ For the rod, as well as for a rigid
body mentioned above,
we consider the left-invariant  metrics, and hence, left cosets 
$S^1\setminus SO(3)$.
For the group of diffeomorphisms we study the right-invariant  metrics
and the space of right cosets $\mathrm{Diff}(S^1)/\mathrm{Rot}(S^1)$.}
(We can see that the only parameter of the rod is its length. In terms
of the inertia operator, this corresponds to the choice of the nonzero
eigenvalue.)

Note that the inertia operator $A=\mathrm{diag}(\lambda, \mu, 0)$ with
$\lambda\not= \mu$ does not correspond to any physical object.
The corresponding degenerate metric on $SO(3)$ is not $S^1$-invariant,
and hence it does not descend to $S^2$.
\smallskip

b) (The Hunter--Saxton equation.) Consider the group of diffeomorphisms
$\mathrm{Diff}(S^1)$ and its quotient
$\mathrm{Diff}(S^1)/\mathrm{Rot}(S^1)$ by the subgroup
of rotations $\mathrm{Rot}(S^1)$.
Consider the degenerate quadratic form on the corresponding Lie algebra
$vect(S^1)$ given by the homogeneous $\dot{H}^1$-energy:
$$
E_0(v\partial_x)=\frac 12\int_{S^1}v_x^2\, dx.
$$
(Similarly,  for  the
Virasoro algebra, consider the energy
$E_0(v\partial_x,b)=\frac 12(\int_{S^1}v_x^2\, dx +b^2).$)

This energy vanishes on constant vector fields. Those fields generate the
subgroup  ${\mathrm{Rot}}(S^1)$ of rotations of the circle $S^1$.
One can see that the form $E_0(v\partial_x)$ is bi-invariant
with respect to the circle action, since the energy is invariant with
respect to translations $x\mapsto x+\mathrm{const}$. (The same holds for the
Gelfand--Fuchs cocycle, and the  energy
$E_0(v\partial_x,b)$ on the extended algebra.)
Hence the corresponding
right-invariant metric on $\mathrm{Diff}(S^1)$ descends to the quotient
$\mathrm{Diff}(S^1)/\mathrm{Rot}(S^1)$. We will see that the geodesics with
respect to this metric are described by the HS equation.

\begin{rem} One hopes that a certain modification of this approach can be
applied to generalized flows in \cite{bre}, where  the fluid particles in 3D
can move freely and independently along one coordinate.
The corresponding subgroup $K$ here might be that of fiberwise
diffeomorphisms along a coordinate.
\end{rem}

\medskip

{\bf Proof of Theorem \ref{coset}.}
First of all we note that the quadratic form $E$ on
$\mathfrak g$ induced from a nondegenerate form on the 
quotient $\mathfrak{g}/\mathfrak{k}$ 
is degenerate exactly along $\mathfrak k$ (i.e., $\mathfrak k$
is its null subspace).

Let $\mathfrak k$ be a Lie subalgebra and  consider the quotient of the
corresponding  groups $G/K$. Note that the restriction of the
energy form to the subgroup $K$ is zero. Indeed, the latter is nothing but
the right translation from the identity of the  energy form on
the subalgebra $\mathfrak k$.

We would like to compare
the right-invariant metric $E_G$ on $G$ at two different points
$k_1 g$ and $k_2 g$ of the same coset $Kg$. Then the element
$\bar k=k_2k_1^{-1}\in K$
sends $k_1$ to $k_2$ by means of the left translation.
This translation also identifies the tangent spaces to $G$ along the
same coset $Kg$, see Fig.2. The energy form $E_G$ is invariant under
this identification, since it is bi-invariant with respect to the action
of elements of $K$. Finally, note that the energy $E_G$ is degenerate
along cosets. (Indeed, it vanishes on $K$, the ``identity coset,''
and it is invariant with respect to right translations, which
shuffle the cosets.)

Therefore, the corresponding energy form descends to the coset space.
$\square$
\bigskip

%%%%%%%%%%%%%%%%%%%%%%%%%%%%%%%%%%%%%%%%%%%%%%%%%%%%%%%
%%%%%%%%%%%%%%%%%%%%%%%%%%%%%%%%%%%%%%%%%%%%%%%%%%%%
\begin{figure}[htb]
\centerline{\epsfbox{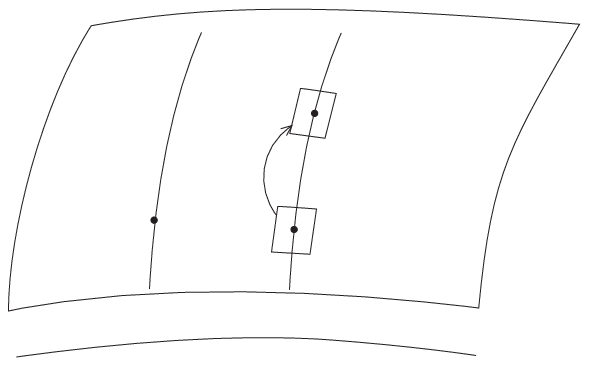}}
\begin{picture} (0, 34)
\put(146,118){$G$} \put(115,33){$G/K$} \put(185,123){$K$}
\put(220,125){$Kg$} \put(170,75){$id$} \put(206,87){$\bar k$}
\put(240,105){$k_2g$} \put(235,70){$k_1g$}
\end{picture}
\vskip -13mm \caption{\small \footnotesize
Defining a right-invariant form on the space
of right cosets $G/K$.} \label{fig4}
\end{figure}
%%%%%%%%%%%%%%%%%%%%%%%%%%%%%%%%%%%%%%%%%%%%%%%%%%%%%%%
%%%%%%%%%%%%%%%%%%%%%%%%%%%%%%%%%%%%%%%%%%%%%%%%%%

\begin{rem} From the Hamiltonian point of view, the 
geodesic picture on a homogeneous space  corresponds to a
Hamiltonian reduction of the non-degenerate case with respect to the
subgroup $K$ action.

More precisely,
consider the cotangent bundle $T^*(G/K)$ of the metric space
$G/K$. 
Similarly to the non-degenerate case, look at the fiber 
$(\mathfrak{g/k})^*$ over the ``identity coset''  $K$. 
This space $(\mathfrak{g/k})^*$ can be naturally identified with the
image $L={\rm im} A\subset \mathfrak g^*$ of the degenerate
inertia operator $A:\mathfrak g\to \mathfrak g^*$, for which  
${\mathfrak k}={\rm ker} A$, (or, equivalently, with the
subspace $L=\mathfrak{k}^\perp\subset \mathfrak{g}^*$,
the annihilator of
$\mathfrak{k}$ in $\mathfrak{g}^*$).

The subgroup $K$, being a
stabilizer of the ``identity coset''and hence of
$L=(\mathfrak{g/k})^*$, 
acts on $L$ respecting the Poisson structure.
Therefore, there is a natural Poisson structure 
on the quotient $L/K:={\rm im}\, A/{\rm ad}^* K$. 

Furthermore, one can define the corresponding Hamiltonian function on $L$
by the same formula as above:
\begin{equation}\label{homHam}
H_L(m)=\frac 12\langle m, A^{-1}m\rangle.
\end{equation}
This function is ${\rm ad}^* K$-invariant and hence is well-defined 
on the quotient $L/K$.

Suppose $A$ generates a right-invariant and ${\rm ad}\, K$-invariant 
metric on the group $G$
(i.e., satisfies the conditions of Theorem \ref{coset}), so that
it makes sense to consider geodesics of the corresponding metric on
$G/K$. 
 
Then the above consideration provides the following  
limiting (degenerate) case of Arnold's theorem
(cf. Definition \ref{linEuler}).
\end{rem}

\begin{thm}\label{homArnold}
The Euler equation, which corresponds to the inertia operator $A$ and
describes the geodesic flow 
on the space $G/K$ of right cosets, has the following Hamiltonian 
form on $L/K = {\rm im}\, A/{\rm ad}^* K$: it is the quotient with 
respect to the $K$-action of the   restriction to 
$L\subset \mathfrak g^*$ of the following
Hamiltonian equation on $\mathfrak g^*$ 
$$
\frac{dm}{dt}=-{\rm ad}^*_{A^{-1}m}m
$$
for $m\in L={\rm im}\, A$.
\end{thm}

Now we are ready to complete the argument showing
the Eulerian nature of the Hunter-Saxton equation.

\medskip

\begin{thm}{\bf (=\ref{homEuler}$'$)}\label{homEuler'} 
The HS equation 
$$%\begin{equation}
v_{txx}=-2 v_x  v_{xx} -v v_{xxx}
$$%\end{equation}
is a well-defined equation on the equivalence classes of
periodic  functions 
$$
\big\{ v(x)\sim v(x+p)+q {\rm ~~for~ any~~ } p,~q \big\}.
$$ 
It describes the geodesic flow on the homogeneous space
$Vir/\mathrm{Rot}(S^1)$ of the Virasoro group modulo rotations
with respect to the right-invariant homogeneous
$\dot{H}^1$ metric.
\end{thm}

{\bf Proof of Theorem \ref{homEuler'}.}
The homogeneous $\dot{H}^1$-metric on the Virasoro algebra is related
to the degenerate  inertia operator $A: vir\to vir^*$ sending
$$
(v\partial_x, b)\mapsto (-(\partial_x^2 v)(dx)^2, b).
$$
Its image $L\subset vir^*$ consists of pairs $(u(dx)^2,a)$,  
where functions $u$ have 
zero mean:
$$
L = \left\{  
(u(dx)^2,a)~ \Big\vert~\int_{S^1} u(x)~dx =0 
\right\}. 
$$
The action of the 
subgroup $K={\rm Rot}(S^1)$ identifies those functions $u$ that
differ by a rotation:
$u(x)\sim u(x+p)$. Thus we come to a Hamiltonian equation on $L/K$.

For explicit calculations, recall that  the present case
corresponds to  setting
$\alpha=0$, $\beta=-1$ in the proof of Theorem \ref{Euler'}.
Choosing these values in Equation (\ref{abEuler}), one arrives
at the HS equation
$$
v_{txx}=-2 v_x  v_{xx} -v v_{xxx}-bv_{xxx}
$$
on $vir$. In order to obtain equation describing the evolution 
of $v$ rather than that of $v_{xx}$ we observe that 
$$
2 v_x  v_{xx} +v v_{xxx}+bv_{xxx}
= 
\Big(vv_{xx}+\frac 12(v_x)^2
+bv_{xx}\Big)_x 
$$
and hence, integrating both sides in $x$, we obtain
$$
v_{tx}=-vv_{xx}-\frac 12(v_x)^2 -bv_{xx}+r,
$$
where $r$ is an arbitrary constant. This constant is uniquely 
determined by the condition that the right-hand-side is a complete 
derivative, i.e., by $\int_{S^1} ((v_x)^2/2 +r_0)dx=0$. Then
$$
v_{t}=-vv_{x}+\partial^{-1}_{x}((v_x)^2/2 +r_0) -bv_{x}+q,
$$
where $q$ is an arbitrary constant.

Whence $v_{t}$ (and hence the evolution of $v$) is defined only up to
addition of a parameter $q$ and a multiple of $v_{x}$. 
(The latter is the velocity of the rotation subgroup: $v_x=\frac{d}{dp}|_{p=0}v(x+p)$.)
This manifests the fact that the evolution of $v$ is defined 
on the equivalence classes $\{v(x+p)+q\}$. 
The inertia operator $A$ sends this
equation on classes to the Hamiltonian equation on the 
quotient $L/K$.

Note that  the equivalence classes above 
 absorb the term $bv_x$ (respectively,
$bv_{xxx}$ for $v_{txx}$: consider a shift $v\mapsto v+{\rm const}$) 
and one obtains the HS equation 
in its standard form (\ref{eq:HS}).~~$\square$

\medskip

One should notice that setting  $b=0$ corresponds to the Euler equation
on the quotient $\mathrm{Diff}(S^1)/\mathrm{Rot}(S^1)$.
Thus in the homogeneous case the consideration of the central extension
does not give anything new, since the Euler equations for
$\mathrm{Diff}(S^1)/\mathrm{Rot}(S^1)$ and ${Vir}/\mathrm{Rot}(S^1)$
are equivalent.

%%%%%%%%%%%%%%%%%%%%%%%%%%%%%%%%%%%%%%%%%%%%%%%%%%%%%%%%

\section{Bihamiltonian structures for the equations}
To formulate our next result we need to recall some generalities
on bihamiltonian systems.

\begin{defn}
Assume that a manifold $M$ is equipped with two Poisson structures
$\{.,.\}_0$ and $\{.,.\}_1$. They are said to {\rm  be compatible}
(or, form a {\rm Poisson pair}) if all of their linear combinations
$\{.,.\}_0+\lambda\{.,.\}_1$ are also Poisson structures.

A dynamical system $dm/dt=F(m)$ on $M$ is called {\rm bi-Hamiltonian} if
the vector field $F$ is Hamiltonian with respect to both structures
$\{.,.\}_0$ and $\{.,.\}_1$.
\end{defn}

Consider the dual space $\mathfrak g^*$ to a Lie algebra $\mathfrak g$.
As we discussed above, it is equipped with the {\it Lie-Poisson structure}:
$$
\{f,g\}_{LP}(m):=\langle [df,dg], m\rangle
$$
where $m\in \mathfrak g^*$ and   $f,g$ are two arbitrary functions on $\mathfrak g^*$.

Now fix a point $m_0$ in  $\mathfrak g^*$.
One can associate to this point another Poisson bracket on $\mathfrak g^*$
as follows, see \cite{ma}.

\begin{defn}
The {\rm constant Poisson bracket} associated
to a point $m_0\in \mathfrak g^*$ is the bracket $\{.,.\}_0$ on the dual
space $\mathfrak g^*$  defined by
$$
\{f,g\}_0(m):=\langle [df,dg], m_0\rangle
$$
for any two smooth functions $f,g$ on the dual space, and any $m\in \mathfrak g^*$.
The differentials $df,dg$ of the functions $f,g$ are taken, as above,
at the point $m$ and are regarded as elements of the Lie algebra itself.
\end{defn}

The constant bracket depends on the choice of the ``freezing''
point $m_0$, while the
Lie--Poisson bracket is defined by the Lie algebra structure  only.
Note that the brackets $\{.,.\}_{LP}$ and $\{.,.\}_0$ coincide
at the point $m_0$ itself, and, moreover, the bivector defining
the constant bracket $\{.,.\}_0$ is the same at all  points $m$.

\begin{prop}
The brackets $\{.,.\}_{LP}$ and $\{.,.\}_0$ are compatible for every
``freezing'' point $m_0$.
\end{prop}

{\bf Proof.}
Indeed, any linear combination $\{.,.\}_\lambda:=\{.,.\}_{LP}+\lambda
\{.,.\}_0$ is again a Poisson bracket, since it is just the
linear Lie--Poisson structure
$\{.,.\}_{LP}$ translated from the origin to the point
$-\lambda m_0$.~~$\square$
\smallskip

\begin{rem}\label{rem:hamconst}
 Explicitly, the Hamiltonian equation on $\mathfrak{g}^*$ with
the Hamiltonian function $f$ and computed with
respect to the constant Poisson structure frozen
at a point $m_0\in \mathfrak{g}^*$ has the following form:
\begin{equation}\label{eq:hamconst}
\frac{dm}{dt}={\rm ad}^*_{df}m_0,
\end{equation}
as a modification of Proposition \ref{hamfield} shows.
\end{rem}

\medskip

Now we can formulate another main result.

\begin{thm}\label{thm:2-h} The Euler equation  (\ref{abEuler})
for the $H^1_{\alpha, \beta}$
-metric (with $\alpha\not=0$) on the  Virasoro group
is bihamiltonian on the  dual  $vir^*$ of the Virasoro algebra.
The corresponding ``freezing'' 
point  in ${vir}^*$ is $(\frac{\alpha}{2}(dx)^2, \beta)$.
\end{thm}

In Appendix we show that the dual space $vir^*=\{ (u(x)(dx)^2, a)\}$
can be thought of as the space of Hill's operators 
$\{ -a\partial_x^2+u(x)\}$.
In these terms 
the above theorem can be stated as follows:
The  $H^1_{\alpha, \beta}$ -metric on $vir$ given by the inertia operator
$\Lambda=\alpha- \beta\partial_x^2$ is bihamiltonian on  $vir^*$ 
with ``freezing'' at the point $\alpha/2-\beta\partial_x^2$.

\medskip

\begin{rem}\label{rem:2-hh} The KdV and CH equations are bihamiltonian
on the  Virasoro dual. The corresponding ``freezing'' 
points $(u_0(dx)^2, a_0)$ in ${vir}^*$ are $(u_0=1/2, a_0=1)$
for the CH equation and
$(u_0=1/2, a_0=0)$ for the KdV equation, see Fig.3, as they are related to
the $H^1$- and  $L^2$-energies, respectively. 

To describe bihamiltonian nature of the HS equation 
on the reduced space $L/K$, 
discussed in Theorems \ref{homArnold}--\ref{homEuler'}, one should
consider the following analog of the constant Poisson structure.
Take the Lie--Poisson structure on  ${vir}^*$ ``frozen'' at the point
$( u_0=0, a_0=1 )$ and then push it forward  to 
the corresponding quotient space
for the HS equation. An alternative way to show integrability 
(rather than the bihamiltonian property) of this equation is to use 
integrability of CH and an infinite-dimensional version of the 
formalism developed in \cite{th}, \cite{gs}. 
\end{rem}

\begin{que}
Which metrics 
on the Virasoro group (or which quadratic forms
on the Virasoro 
 algebra) correspond to the bihamiltonian system
on  ${vir}^*$ with ``freezing'' at a point $( u_0(x)(dx)^2, a_0 )$
for non-constant $u_0(x)$?\footnote{After this paper was submitted,
I.~Zakharevich  found 
a  formula for the corresponding quadratic form on the Lie algebra in terms of
Bloch solutions of the operator  $ -\partial_{x}^2+u(x)$, \cite{za}.}
For which $u_0(x)$  are these metrics 
positive definite
(i.e., Riemannian rather than pseudo-Riemannian)?
\end{que}

%%%%%%%%%%%%%%%%%%%%%%%%%%%%%%%%%%
%%%%%%%%%%%%%%%%%%%%%%%%%%%%%%%%%%
\vskip 2mm
\begin{figure}[htb] \centerline{\epsfbox{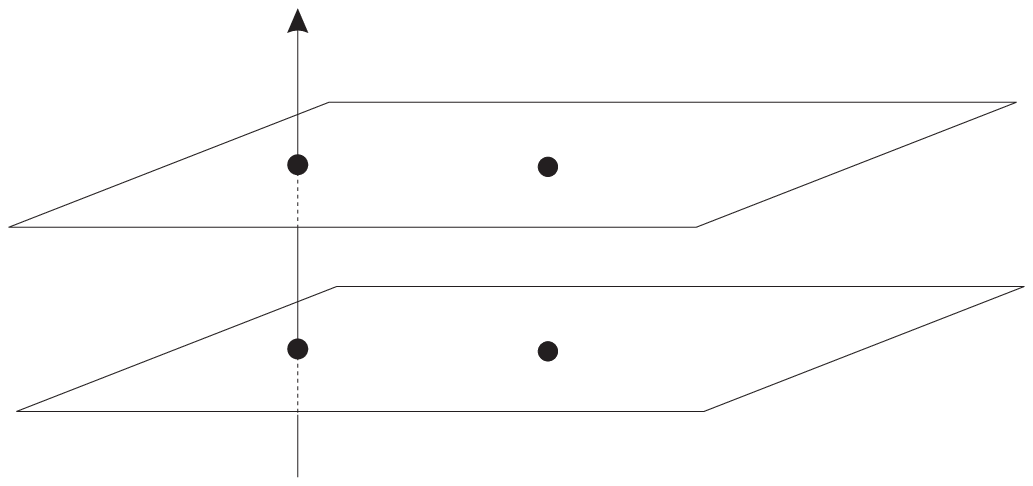}}
\begin{picture} (0,0)
\put(240,102){CH} \put(330,82){$\{u(x)\}$}
\put(200,142){$vir^*=\{-a\partial_{x}^2+u(x)\}$} \put(146,142){$a$}
\put(150,42){$0$} \put(150,96){$1$} \put(170,102){HS}
\put(240,48){KdV}
\end{picture}
 \vskip -5mm \caption{\small \footnotesize
Locations
of the ``freezing'' points for the KdV, CH, and HS equations in
the Virasoro \break dual $vir^*$.} \label{fig3}
\end{figure}
%%%%%%%%%%%%%%%%%%%%%%%%%%%%%%%%%%%%
%%%%%%%%%%%%%%%%%%%%%%%%%%%%%%%%%%%%

{\bf{Proof of Theorem \ref{thm:2-h}}.}
Let $F(u,a)$ be a function on ${vir}^*$ and  let
$
(v\partial_x , b) = ({\delta F}/{\delta u}, {\delta F}/{\delta a})
$
be a (variational) derivative of $F$ at $(u(dx)^2,a)$.
Then the Hamiltonian equation with Hamiltonian function $F$ computed with
respect to the constant Poisson structure ``frozen'' at
$(u_0(dx)^2 , a_0)$
has the form
$$
\frac{d}{dt}(u(dx)^2, a)
=
{\rm ad}^{*}_{(v\partial_x,b)}(u_0(dx)^2 , a_0)
=
\left( (2u_0 v_x - a_0 v_{xxx})(dx)^2 , 0\right).
$$
(Here we used Remark \ref{rem:hamconst}, the explicit form
(\ref{virad*}) of ${\rm ad}^*$ for
the Virasoro algebra,  and the fact that $u_0 = \text{const}$.)

For the first component of this equation one has
$$
\frac{du}{dt}
=
( 2u_0 - a_0\partial_x^{2})\,\partial_x\left(\frac{\delta F}{\delta u}\right).
$$
Setting $u_0 = \alpha/2$ and $a_0 = \beta$
we obtain
$ 2u_0 - a_0\partial_x^{2}
=
\alpha - \beta\partial_x^{2}
=
\Lambda,$
and this simplifies the equation to
\begin{equation}\label{eq:Lam}
\frac{du}{dt}
=
\Lambda \left( \partial_x\left(\frac{\delta F}{\delta u}\right)\right).
\end{equation}
To prove the Theorem one needs
to show that  for any 
$\alpha\not =0$ and any  $ \beta$
the Euler equation (\ref{abEuler}) can also be expressed 
in the form (\ref{eq:Lam})
for an appropriate Hamiltonian function $F$.

Next, consider the Hamiltonian function $F$ of the form
$$
F(u,a)
=
- \int \left(
\frac{\alpha}{3} (\Lambda^{-1} u)^3
+
\frac{1}{4} (\Lambda^{-1}u)^{2}u
+
\frac{a}{2} (\Lambda^{-1}u_{x})^2
\right) \, dx.
$$
(The operator $\Lambda:=\alpha - \beta\partial_x^{2}$ is invertible for
$\alpha\not =0$.)
By definition, the variational derivative
$\left({\delta F}/{\delta u}, {\delta F}/{\delta a}\right)\in vir $
 of the functional $F$ is determined by the following identity satisfied
for any
$ (\xi(dx)^2, c)\in vir^*$:
$$
\left\langle (\xi(dx)^2, c) ,
\left(\frac{\delta F}{\delta u},\frac{\delta F}{\delta a}\right)
\right\rangle
=
\frac{d}{d\epsilon}\Big\vert_{\epsilon =0} F(u + \epsilon \xi ,a+\epsilon c).
$$
Since we need only the partial variational derivative
${\delta F}/{\delta u} $, we compute:
$$
\frac{d}{d\epsilon}\Big\vert_{\epsilon =0} F(u + \epsilon \xi ,a)=
$$
$$
\frac{d}{d\epsilon}\Big\vert_{\epsilon =0}
\int
\left(
\frac{\alpha}{3} (\Lambda^{-1} (u+ \epsilon \xi))^3
+
\frac{1}{4} (\Lambda^{-1}(u+ \epsilon \xi))^{2}(u+ \epsilon \xi)
+
\frac{a}{2} (\Lambda^{-1}(u+ \epsilon \xi)_{x})^2
\right) \, dx
=
$$
$$
-  \int \xi\cdot
\left(
\alpha \Lambda^{-1} (\Lambda^{-1}u)^2
+
\frac{1}{4} (\Lambda^{-1}u)^{2}
+
\frac{1}{2} \Lambda^{-1}( (\Lambda^{-1}u) u)
-
a \Lambda^{-2}u_{xx}
  \right)\,dx.
$$
Thus, we have found that
$$
\frac{\delta F}{\delta u}=-\left(
\alpha \Lambda^{-1} (\Lambda^{-1}u)^2
+
\frac{1}{4} (\Lambda^{-1}u)^{2}
+
\frac{1}{2} \Lambda^{-1}( (\Lambda^{-1}u) u)
-
a \Lambda^{-2}u_{xx}
  \right)\partial_x
$$
Now, we substitute the variational derivative ${\delta F}/{\delta u}$
into Equation (\ref{eq:Lam}) and then rewrite the obtained equation
on the algebra $vir$, rather than on its dual $vir^*$.
The latter corresponds to rewriting the equation in terms of
$(v\partial_x , b)=A^{-1}(u(dx)^2,a)$, i.e., in terms of
an unknown function $v=\Lambda^{-1} u$ and setting $b=a$.

Finally, applying $\Lambda$ to the equation we obtain
$$
\Lambda\left(\frac{dv}{dt}\right)= -\left(2\alpha\, vv_x+\frac 12 \Lambda(vv_x)
+\frac 12 v\Lambda v_x+\frac 12 v_x\Lambda v- bv_{xxx}\right).
$$
Recalling that $\Lambda =\alpha-\beta\partial_x^2$ and collecting the terms
we recover Equation (\ref{abEuler}).~~$\square$

\bigskip

To explain in what sense the above equations are
generic bihamiltonian systems on $vir^*$
(obtained by the freezing argument method), we need to
consider symplectic leaves of the above Poisson structures.
\bigskip

%%%%%%%%%%%%%%%%%%%%%%%%%%%%%%%%%%%%%%%%%%%%%%%%%%%%

\section{Hierarchies of Hamiltonians from compatible \\structures}

Recall that the symplectic leaves, i.e.,
maximal non-degenerate submanifolds, of the Lie--Poisson structure
are the coadjoint orbits of the group action on  $\mathfrak g^*$
(see, e.g., \cite{kir} or this also follows from 
Proposition \ref{hamfield}).
Therefore the functions  constant on symplectic
leaves (called {\it Casimir functions})
of the Lie--Poisson bracket are those functions
on the dual space $\mathfrak g^*$ that are invariant under the coadjoint action.
The tangent plane to the group coadjoint orbit at the point $m_0$, as
well as  all  the planes in  $\mathfrak g^*$ parallel to this tangent plane
are the symplectic leaves of the constant bracket frozen at the point $m_0$ .

\begin{defn} The codimension of the coadjoint orbit passing through $m_0$
will be called the {\rm codimension of the Poisson pair} $\{.,.\}_0$
and $\{.,.\}_{LP}$.
\end{defn}

It turns out that there are no Poisson pairs of codimension 0 or 1 in 
the (smooth) Virasoro dual $vir^*$,
and that the Poisson pairs of codimension 2 can all be 
classified.

\begin{thm}\label{pairclass}
All Poisson pairs $\{.,.\}_0$ and $\{.,.\}_{LP}$
on $vir^*$ of codimension 2 belong to one of three classes 
according to the orbit type of the  ``freezing'' point $(u_0(dx)^2, a_0)$.
These classes can be  represented by the points  a) $((dx)^2/2,1)$, 
 b) $((dx)^2/2,0)$, and c) $(0,1)$. 
\end{thm}

{\bf{Proof.}} First we observe that for  any Lie algebra $\mathfrak g$
the list of coadjoint orbits in $\mathfrak g^*$
provides  the list of normal forms for the constant and 
Lie--Poisson pairs as well. Indeed, let a point  $\bar m_0\in \mathfrak g^*$ 
be a normal form for all points $m_0$ belonging to the same coadjoint orbit
as $\bar m_0$.
This means that $m_0$ can be mapped by
the group coadjoint action to its normal form:
$\bar m_0={\rm Ad}^*_g m_0$. The group action
${\rm Ad}^*_g: \mathfrak g^*\to \mathfrak g^*$ 
is a linear operator on $\mathfrak g^*$,  which 
preserves the Lie--Poisson bracket on ${\mathfrak g}^*$.
It also maps  the constant bracket frozen at $m_0$ to
that frozen at $\bar m_0$. Thus the group action ${\rm Ad}^*_g$ 
sends one Poisson pair to the other.

The proof  of Theorem \ref{pairclass}\, is based
on the Virasoro orbit classification, which we recall in  Appendix.
Notice that a
``cocentral value'' $a_0$ is invariant on the orbits in $vir^*$.
This allows us to fix  $a_0$ and consider
the orbits in the  hyperplane  $\{(u(x)(dx)^2,a)|~a=a_0\}$.

As shown in Corollary \ref{normforms} of Appendix, 
for $a_0\not=0$ there are  exactly two types of Virasoro orbits
of codimension 1 in this hyperplane that correspond to
either a generic or a  Jordan $2\times 2$ block holonomy for 
Hill's operator $-a_0\partial^2_x+u(x)$.
If we discard the discrete invariant of these orbits (see Appendix),
representatives for
their normal forms can be chosen  as stated in the a) and b) parts of the
Theorem.

If $a_0=0$ there is just one orbit type
of codimension 1 in the corresponding hyperplane, and represented by c),
see Remark \ref{kdvfreez}.
Note that in the whole dual space those orbits have codimension 2, taking 
extra dimension $a$ into account.~~$\,\square$

\medskip

Given a Poisson pair one can generate a bihamiltonian dynamical system,
by producing a sequence of  Hamiltonians in involution, according to the
following {\it Lenard scheme}.
Let $\{.,.\}_\lambda:=\{.,.\}_0+\lambda\{.,.\}_1$
be the Poisson bracket on a manifold $M$
for any $\lambda$. Denote by $h_\lambda$ its Casimir function on $M$
parameterized by $\lambda$.
This means that $\{h_\lambda,f\}_\lambda=0$ for any function $f$.
Expand $h_\lambda$ in  a power series: $h_\lambda=h_0+h_1\lambda+...$,
where each coefficient $h_j$ is a function on $M$. The following 
theorem is well-known.

\begin{thm}\label{biham}
The functions $h_j,~j=1,2\dots$ are Hamiltonians of a hierarchy of
bihamiltonian systems.
In other words, each function $h_j$ generates the Hamiltonian field
$X_j$  with respect to the
Poisson bracket $\{.,.\}_1$ (i.e., $X_j$  satisfies
$L_{X_j}f= \{h_j,f\}_1 $ for any $ f$), which is also Hamiltonian for the
other bracket
$\{.,.\}_0$ with Hamiltonian function $-h_{j+1}$ (i.e.,
$L_{X_j}f= -\{h_{j+1},f\}_0 $ for any $ f$).
Other functions $h_i,~i\not=j$ are first integrals of
the corresponding  dynamical systems $X_j$.
\end{thm}

In other words, the  functions
$h_j,~j=0,1\dots$ are in involution with respect to each of the
two Poisson brackets $\{.,.\}_0$ and $\{.,.\}_1$.

{\bf Proof.}
Substituting the power series for $h_\lambda$
 into the Casimir condition  we obtain:
$$
0=\{h_\lambda,f\}_\lambda=\{h_0+h_1\lambda+...,f\}_0
+\lambda\{h_0+h_1\lambda+...,f\}_1.
$$
Collecting the terms at $\lambda^0,~\lambda^1,~\lambda^2,\dots$ we obtain
a sequence of identities:
$$
\{h_0,f\}_0=0,~~\{h_1,f\}_0+\{h_0,f\}_1=0,~~\{h_2,f\}_0+
\{h_1,f\}_1=0,\dots
$$
for any function $f$.  The first identity expresses the fact that $h_0$
is a Casimir function for the bracket $\{.,.\}_0$.
The next one says that the Hamiltonian field for $h_1$ with respect to
$\{.,.\}_0$ coincides with the Hamiltonian field for $-h_0$ and the bracket
$\{.,.\}_1$, and so on.

To see that every function $h_i$ is a first integral for the equation
generated by $h_j$ with respect to each bracket, we check that
$\{h_i,h_j\}_k=0, ~k=0,1.$ Indeed, e.g., if $i<j$
$$
\{h_i,h_j\}_1=-\{h_{i+1},h_j\}_0=\{h_{i+1},h_{j-1}\}_1=\dots=0,
$$
since we finally obtain the bracket (either $\{.,.\}_0$ or $\{.,.\}_1$)
of one of the functions $h_l$ with itself.
$\,\square$

%\bigskip

%\bigskip

Thus the choice of Casimir functions $h_\lambda$ determines the
corresponding (hierarchy of) dynamical systems. 
By combining Theorems
\ref{pairclass} and \ref{biham} we get the following

\begin{cor}\label{classpairs}
The three types of Poisson pairs of codimension 2 on the Virasoro algebra
correspond to the three integrable systems: CH, KdV, and HS.
These three systems represent all generic Hamiltonian systems
on $vir^*$ (modulo the ambiguity in the choice of Casimir),
which can be integrated by the freezing argument method.
\end{cor}

\begin{rem}\label{ell-hyp}
The corresponding ``freezing'' points in $vir^*$ represent all three types of
the Virasoro coadjoint orbits of codimension 2.

If the ``freezing'' point $(u_0(dx)^2, a_0)$ is generic, 
one obtains an equation ``equivalent'' to
the CH equation. In this sense, the CH equation is the most general equation,
which is encountered by applying the freezing argument method of integration;
this is the case a) in Theorem \ref{pairclass}.

Two other equations can be recovered by confining the ``freezing'' point
to special hypersurfaces in $vir^*$.
(In turn, these hypersurfaces are foliated
into coadjoint orbits. Those orbits are of codimension 1 in the hypersurfaces,
and hence of total codimension 2 in $vir^*$. The classification of
the orbits will be discussed in detail in Appendix.)
A generic point on the hyperplane $a_0=0$ produces the KdV equation;
see the case b) in Theorem \ref{pairclass}.
The case c) in the same theorem corresponds to the HS equation
if we consider  the cone-like Virasoro orbits in the $a_0\not=0$-case
(see Appendix and Fig.4).
The latter  (e.g., ``freezing'' at the point $(u_0(dx)^2, a_0)=(0,1)$)
corresponds to the Euler equation with a
degenerate metric on the group, which we discussed in Section  \ref{homsect}.

\smallskip
One could also consider a more subtle Virasoro orbit 
classification, where
one distinguishes between the two types of generic orbits in  $vir^*$: 
hyperbolic and elliptic ones, according to the eigenvalues of the monodromy, as
well as between the orbits which differ by the discrete invariant
(see Corollary \ref{normforms}). 
While all elliptic orbits with arbitrary values of the
discrete invariant can be represented by Hill's operators with {\it constant 
coefficients}, just one of the hyperbolic and one of the Jordan block classes
has such representatives, while others do not.

Note that the bihamiltonian equations corresponding to  elliptic 
orbits with different discrete invariants are almost the same:
these are the CH equations with different coefficients. 
It would be very interesting to see whether an analogous similarity holds
for the hyperbolic and  Jordan block orbits with different discrete invariants.
\end{rem}

\bigskip

\begin{rem}
When the symplectic leaves of the bracket
$\{.,.\}_\lambda$ are of codimension 1, then the choice of a Casimir function
is essentially unique for every $\lambda$.
(Any two Casimir functions for every fixed $\lambda$ are
functionally dependent.) 
Therefore, the choice of the Poisson pair itself
defines the bihamiltonian system (modulo the mentioned functional
dependence of the initial Hamiltonian),
provided that the symplectic leaves are hypersurfaces.

This is indeed the case for the Virasoro coadjoint orbits
discussed above,
which have codimension 1 for a fixed cocentral value
$a$. It turns out that a natural Casimir function $h_\lambda(u(x)(dx)^2,a)$
corresponding to the KdV Poisson pair
on the dual space to the Virasoro algebra is the trace of the
monodromy operator associated to Hill's operator
$-a\partial_x^2+u(x)-\lambda^2$.
It generates the first integrals 
of the KdV equation, see Remark \ref{kdvint} in Appendix. 
Similarly, one can expand  Casimir functions for the two other
integrable cases, the CH and HS equations.

Note that for orbits of higher codimension one can start with several
Casimirs and consider several Lenard 
schemes to generate the sequences of Hamiltonians. 
\end{rem}

\bigskip

\begin{rem}
A generic Virasoro coadjoint orbit ${\rm Diff}(S^1)/S^1$ can be 
equipped with a complex structure and a two-parameter family of 
compatible (pseudo-) Kahler metrics \cite{kir2}. 

This  family of Kahler 
metrics has a simple origin: a generic Virasoro orbit has 
codimension 2, i.e., it is locally included in a two-parameter
family of orbits, each equipped with its own symplectic structure
compatible with the complex structure.
Alternatively, one could consider a two-parameter 
family of symplectic structures on the same orbit, 
given by the Hamiltonian operators $a\partial_x^3+b\partial_x$.

It turns out that the restriction of  the two-parameter family of
$H^1_{\alpha, \beta}$ -metrics on $vir^*$ to a generic Virasoro orbit
 ${\rm Diff}(S^1)/S^1$ coincides with the family of Kahler metrics
on the orbits. Proof is achieved by comparison with  formula (7) in \cite{kir2}
for those homogeneous metrics at one point of the orbit.

This is yet another fact manifesting a special role of the 
$H^1_{\alpha, \beta}$ -metrics in Virasoro geometry.
\end{rem}

\bigskip

%%%%%%%%%%%%%%%%%%%%%%%%%%%%%%%%%%%%%%%%%%%
%%%%%%%%%%%%%%%%%%%%%%%%%%%%%%%%%%%%%%%%%%%

\section{Appendix: Classification of Virasoro  orbits}

In this section we recall  the classification result for Virasoro
coadjoint orbits (see, e.g., \cite{kir}, \cite{seg} or the book
\cite{gr}). The dual spaces to
the infinite-dimensional Lie algebras considered below are
always understood as smooth duals, i.e. identified with
appropriate spaces of smooth functions.

\medskip

{\it A. Classification of quadratic differentials.}
We start with the non-extended group of  diffeomorphisms
of the circle. Let $\mathrm{Diff}(S^1)$ be
the group of all orientation-preserving diffeomorphisms of $S^1$
and let $vect(S^1)$ be its Lie algebra.

\begin{prop}\label{quaddif} {\rm{\cite{kir}}}
The dual space $vect(S^1)^*$ is naturally identified with the space of
quadratic differentials $\{u(x)(dx)^2\}$ on the circle. 
The pairing is given by the formula:
$$
\langle u(x)(dx)^2,~v(x)\partial_x) \rangle=\int_{S^1} u(x)v(x)~dx
$$
for any vector field $v(x)\partial_x\in vect(S^1)$.
The coadjoint action
coincides with the action of a diffeomorphism on the quadratic differential:
for a diffeomorphism $\varphi\in \mathrm{Diff}(S)$ the action is
$$
{\rm Ad}^*_\varphi:~~ u\,(dx)^2 \mapsto u(\varphi)\cdot\varphi_x^2\,(dx)^2=
u(\varphi)\cdot(d\varphi)^{2}.
$$
\end{prop}

Hence, for instance, if $u(x)>0$ for all $x\in S^1$ then 
the square root $\sqrt{u(x)(dx)^2}$  transforms under
a diffeomorphism as a differential 1-form. In particular,
$\Phi(u(x)(dx)^2)=\int_{S^1} \sqrt{u(x)}\,dx$ is a Casimir function 
(i.e., an 
invariant of the coadjoint action). One can see that there is only one
Casimir function in this case, since  the corresponding
orbit has codimension 1
in the dual space $vect(S^1)^*$.
Indeed, a diffeomorphism action sends
the quadratic differential
${u(x)}(dx)^2$ 
to the constant quadratic
differential $C(dx)^2$, where the constant $C$ is the average value
of the 1-form $\sqrt{u(x)}\,dx$ on the circle:
$$
C=\frac{1}{2\pi}\int_{S^1} \sqrt{u(x)}\,dx.
$$

On the other hand, if a differential ${u(x)}(dx)^2$
changes sign on the circle, then the integral
$\int_a^b \sqrt{|u(x)|}\,dx$, evaluated
between  any two consecutive zeros $a$ and $b$
of the function $u(x)$, is invariant.
In particular, the coadjoint  orbit of such a  differential ${u(x)}(dx)^2$
has necessarily  codimension higher than 1.

\begin{rem}\label{kdvfreez}
In our study of the KdV equation we pick the ``freezing'' point in the dual
space $vect(S^1)^*$ to be $C=1/2$.
(Actually, we consider the dual space to the Virasoro algebra, but we choose
the cocentral term equal to zero, so that the ``freezing'' point
$(u_0(dx)^2, a_0)=((dx)^2/2, 0)$ belongs to
the dual space to the Lie algebra of vector fields.)
Other values of $C$ give equivalent equations, different from the
KdV by scaling only.
\end{rem}

%%%%%%%%%%%%%%%%%%%%%%%%%%%%%%%%%%%%%%%%%%%

\bigskip

{\it B. Virasoro dual and Hill's operators.}
Let  $vir$ be the Virasoro algebra.
We can think of its dual space as the space of pairs
$vir^*=\{(u(x)(dx)^2, a)\}$
consisting of a quadratic differential and a real number (cocentral term).
It is more convenient, however,
to regard such pairs as Hill's operators, i.e. differential operators
$-a\partial_x^2+u(x)$, as we will see below.

\begin{prop}
The Virasoro coadjoint group action is given by the formula
\begin{equation}\label{viract}
{\rm Ad}^*_{(\varphi,\, b)}:~~ (u\,(dx)^2,\, a) \mapsto
\left(u(\varphi)\cdot\varphi_x^2\,(dx)^2-a S(\varphi) \,(dx)^2,\, a\right),
\end{equation}
where
$$
S(\varphi)=\frac{\varphi_x\varphi_{xxx}-
\frac 32\varphi_{xx}^2}{\varphi_x^2}
$$
is the Schwarzian derivative of $\varphi$.
\end{prop}

This group action on Hill's operators:
$$
{\rm Ad}^*_{(\varphi,\, b)}:~~  -a\partial_x^2+u(x)\mapsto
 -a\partial_x^2+u(\varphi)\cdot\varphi_x^2-a S(\varphi)
$$
has the following geometric interpretation
(see, e.g., \cite{seg}, \cite{kir}, \cite{ovs}).
Fix $a=-1$ and consider Hill's operators of the form $\partial_x^2+u(x),
 ~x\in S^1$.
Let $f$ and $g$ be two independent solutions of the corresponding
differential equation
$$
(\partial_x^2+u(x))y=0
$$
for an unknown function $y$.
Although the equation has periodic coefficients,
the solutions need  not necessarily be periodic, but instead
are defined over $\mathbb R$.
Consider the ratio $\eta:=f/g :{\mathbb R}\to {\mathbb RP}^1$.
\smallskip

\begin{prop}
The potential $u$ is (one half of) 
the Schwarzian derivative of the ratio $\eta$:
$$
u=\frac{S(\eta)}{2}.
$$
\end{prop}

{\bf Proof.}
Note that the Wronskian $W(f,g):=fg_x-f_xg$ is constant, since
it should satisfy the differential equation $W_x=0$.
Here we normalize $W$ by setting $W=-1$. This additional condition
allows one to find the potential $u$ from the ratio $\eta$.
Indeed, first one reconstructs the solutions $f, g$
from the ratio $\eta$ by differentiating:
$$
\eta_x=\frac{f_xg-fg_x}{g^2}=\frac{-W}{g^2}=\frac{1}{g^2}\,.
$$
Therefore,
$$
g=\frac{1}{\sqrt{\eta_x}}, \qquad f=g\cdot\eta=\frac{\eta}{\sqrt{\eta_x}}.
$$
Given two solutions $f$ and $g$, one immediately finds the corresponding
differential equation they satisfy
by writing out the following $3\times 3$-determinant:
$$
\det\left[\begin{array}{ccc}
               y&f&g\\
y_x&f_x&g_x\\
y_{xx}&f_{xx}&g_{xx}
               \end{array} \right] = 0
$$
Since  $f$ and $g$ satisfy the equation 
$y_{xx}+u\cdot y=0$, one obtains from the determinant above that
$$
u=-\det\left[\begin{array}{cc}
f_x&g_x\\
f_{xx}&g_{xx}
               \end{array} \right]\,.
$$
The explicit formula for $u$ expressed in terms of $\eta$
turns out to be one half of the Schwarzian derivative of $\eta$.~~$\square$

\medskip

\begin{cor}
The Schwarzian derivative $S(\eta)$
is invariant with respect to a M\"obius transformation
$\eta\to(a\eta+b)/(c\eta+d)$, where $ad-bc=1$.
\end{cor}

{\bf Proof.}
Indeed, for a given potential $u$ the solutions $f$ and $g$
of the corresponding differential equation are not defined uniquely,
but up to a transformation
of the pair $(f,g)$ by  a matrix from $SL_2(\mathbb R)$.
Then the ratio $\eta$ changes by a M\"obius transformation.~~$\square$

\medskip

\begin{prop}
The Virasoro coadjoint action of a diffeomorphism $\varphi$
on the  potential $u(x)$
gives rise to a diffeomorphism change of coordinate in the ratio $\eta$:
$$
\varphi:~~ \eta(x)\to \eta(\varphi(x))
$$
\end{prop}

{\bf Proof.} We look at the corresponding infinitesimal action
on the solutions of the differential equation
$(\partial_x^2+u(x))y=0 $.
For $\varphi(x)=x+\epsilon v(x)$  close to the identity,
consider the  infinitesimal Virasoro action of $\varphi$ 
on the potential $u(x)$:
$$
u\mapsto u+\epsilon\cdot \delta u, ~~{\rm where}~~
\delta u=2uv_x+u_xv-\frac 12 v_{xxx},
$$
(cf. formula (\ref{virad*}) for $a=\frac 12$). It 
is consistent with the following action 
on a solution $y$ of the above differential equation:
$$
y\mapsto y+\epsilon\cdot \delta y, ~~{\rm where}~~
\delta y=-\frac 12 yv_x+y_xv.
$$
The consistency means that 
$(\partial_x^2+u+\epsilon\cdot \delta u)(y+\epsilon\cdot \delta y)
=0+{\cal O}({\epsilon}^2)$.

Note that the action $\epsilon\cdot \delta y
=\epsilon\cdot(-\frac 12 yv_x+y_xv)$ is an infinitesimal version of the
following action of the diffeomorphism $\varphi(x)=x+\epsilon v(x)$ on $y$:
$$
\varphi: ~~y\mapsto y(\varphi)(\varphi_x)^{-1/2}.
$$
Thus solutions to Hill's equation  transform as  forms of degree $-1/2$. 
Therefore the ratio $\eta$ of two solutions
transforms as a function under a diffeomorphism action.~~$\square$
\medskip

In short, 
to calculate the  coadjoint action on the  potential $u$ one can first
pass from this potential to the ratio of two solutions, 
then change the variable
in the ratio, and finally take the Schwarzian derivative of the new ratio to
reconstruct the new potential ${\rm Ad}^*_{(\varphi,\, b)}u$.

\bigskip

All of the above considerations of Hill's operators were local in $x$.
To describe the Virasoro orbits, we now recall that
$u(x)$ is defined on a circle.

\begin{thm} ~
{\rm{\cite{seg}, \cite{kir}}}
The coadjoint Virasoro orbits (for a given cocentral
term $a\not=0$) are enumerated by the conjugacy classes
in $(\widetilde{SL}_2({\mathbb R})\setminus
 \{\text{id}\})/{\mathbb Z}_2$,
the universal covering of $SL_2(\mathbb R)$ without the identity 
and modulo the ${\mathbb Z}_2$-action.
\end{thm}

{\bf Proof.} For a periodic potential $u$
the solutions of $(\partial^2_x+u)y=0$ are quasiperiodic. 
In other words, the boundary  values of the
fundamental set of solutions $F:=(f,g)$ 
 on $[0, 2\pi]$ are related by a holonomy matrix
${M}\in SL_2(\mathbb R)$: $F(2\pi)=F(0)M$.
(Similarly, one can consider the ``projective solution''
$\Theta:=(\eta, \eta_x)$, which consists of the solution ratio $\eta$ and
its derivative $\eta_x$  
 with the correponding holonomy ${\mathcal M}$ now in 
$PSL_2(\mathbb R)$.)
This holonomy matrix $M$ changes to a conjugate matrix 
if $x_0=0$ is replaced by any point $x_0\in S^1$
or if $F$ is replaced by another system of solutions.

Now regard  the ratio $\eta=f/g$ as a map
$\eta:~[0, 2\pi]\to{\mathbb RP}^1$ describing a motion (``rotation'')
along the circle ${\mathbb RP}^1\approx S^1$.
One can see that the condition $\eta_x\not=0$ is equivalent to
the condition $W\not=0$ on the Wronskian.
Choosing  the negative 
sign of the Wronskian, $W<0$, we can assume that
the rotation always goes in the positive direction:
$\eta_x=-W/g^2>0$.
 
By a diffeomorphism change of the parameter
$x\mapsto \varphi(x)$, one can always turn  the map
$\eta:~[0, 2\pi]\to{\mathbb RP}^1$ into a {\it uniform}  rotation along
${\mathbb RP}^1$, while keeping the boundary values of $\eta(x)$
on the segment $[0, 2\pi]$ satisfying
the holonomy relation $\Theta(2\pi)=\Theta(0){\mathcal M}$.
Furthermore, the number of rotations (the ``winding number") for the map
$\eta:~[0, 2\pi]\to{\mathbb RP}^1$ does not change under a reparametrization
by a diffeomorphism $\varphi$.
Thus the orbits of the maps $\eta$ (or, equivalently, of the
potentials $\{u(x)\}$) are described by the conjugacy classes of matrices in
the universal covering of $SL_2(\mathbb R)$. 
The choice  in the sign of the Wronskian reflects 
the ${\mathbb Z}_2$-action on this universal covering.

Finally, note that the identity matrix in the universal covering
 $\widetilde{SL}_2({\mathbb R})$ (or in its projectivization)
cannot be obtained as a holonomy matrix for the maps
$\eta:~[0, 2\pi]\to{\mathbb RP}^1$. Indeed, any map $\eta$
starts at the identity and goes in the positive direction. 
Thus, no matter how slow the rotation, one always moves out 
from the identity.~~$\square$
 
\bigskip

\begin{cor}\label{normforms}
The Virasoro orbits in the hyperplane
$\{-a\partial^2_x+u(x)~| \,\, a=a_0\}\subset vir^*$
with fixed $a\not=0$ are classified by the
Jordan normal form of matrices in $SL_2(\mathbb R)$ and a positive integer
parameter, winding number. 
Matrices in the group
$SL_2(\mathbb R)$ split into three classes, whose normal forms are
the exponents of the following three classes:
$$
i)~~\left[\begin{array}{cc}
\mu&0\\
0&-\mu
               \end{array} \right]\,,~~
ii)~~\left[\begin{array}{cc}
0&\pm 1\\
0&0
               \end{array} \right]\,,~~
{\text{ and } ~~~~}
iii)~~\left[\begin{array}{cc}
0&0\\
0&0
               \end{array} \right]\,
$$
in the corresponding Lie algebra $sl_2(\mathbb R)$, see Fig.4. 
The Virasoro orbit containing
 Hill's operator $\partial_x^2+u(x)$ has the
codimension in the hyperplane $\{a=a_0\}\subset vir^*$
that is equal to the codimension  in $SL_2(\mathbb R)$ of (the conjugacy class
of) the
holonomy matrix $M$ corresponding to this Hill's operator.
This codimension is 1 for the classes $i)$ and $ii)$,
and it is 3 for $iii)$ independ of the integer
parameter.
\end{cor}

\begin{rem}
Note that the class of real matrices whose exponents have
the normal form $i)$ decomposes
into rotation matrices (with $\mu\in i\mathbb R$) and  hyperbolic rotations
(with $\mu\in\mathbb R$), the elliptic and hyperbolic cases, cf. 
Remark \ref{ell-hyp}.
In the paper we regard these cases as belonging to the same general class
(and similarly 
we do not distinguish between the cases 
$\pm 1$ in type $ii)$), 
since we are interested in the algebraic (rather than geometric)
question of constructing the
corresponding integrable equations.

One should also mention that the equality of codimension of the 
Virasoro coadjoint orbits in $vir^*$ and  (the conjugacy class
of) the corresponding holonomy matrices in 
in $SL_2(\mathbb R)$ can be seen by checking the smooth dependence
on a parameter in the above classification.
(The versal deformations of the orbits can be given in terms of the 
Jordan--Arnold normal forms of the holonomy matrices depending on a parameter,
cf. \cite{ok2}.)
Alternatively, one can find the dimension of the corresponding 
stabilizers, see \cite{kir}, \cite{seg}. 

Regarded as homogeneous spaces, the orbits of type $i)$ are often denoted by 
${\rm Diff}(S^1)/S^1$, the notation ${\rm Diff}(S^1)/\mathbb R^1$
stands for $ii)$ (and sometimes for the case $\mu\in i{\mathbb R}$ in $i)$), 
and ${\rm Diff}(S^1)/SL_2(\mathbb R)$ corresponds to $iii)$.
\end{rem}

\medskip

%%%%%%%%%%%%%%%%%%%%%%%%%%%%%%%%%%%%%%%%%%%%%%%%%
%%%%%%%%%%%%%%%%%%%%%%%%%%%%%%%%%%%%%%%%%%%%%%
\begin{figure}[htb]
\centerline{\epsfbox{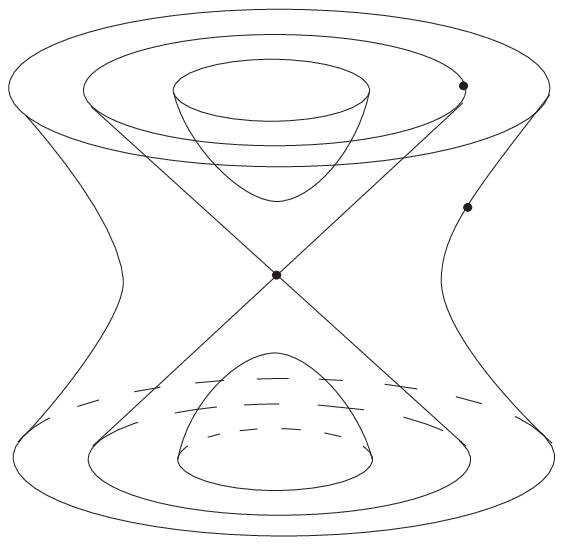}}
\begin{picture} (0, 49)
\put(208,122){$0$} \put(110,120){$sl_2(\mathbb R)$}
\put(283,143){$M_1$} \put(256,176){$M_2$}
\end{picture}
\vskip -16mm \caption{\small \footnotesize
The points $M_1,~M_2$ and 0 in $sl_2(\mathbb R)$
(which is a local picture of $SL_2(\mathbb R)$) correspond to the
Virasoro orbits of the types $i),~ii)$ and $iii)$, respectively.}
\label{fig2}
\end{figure}
%%%%%%%%%%%%%%%%%%%%%%%%%%%%%%%%%%%%%%%%%%%%%%%
%%%%%%%%%%%%%%%%%%%%%%%%%%%%%%%%%%%%%%%%%%%%%%%%%

\begin{rem} For applications to bihamiltonian systems we would like to describe
all  points in $vir^*$ belonging to orbits of codimension at most 2.
As we have shown above, in the smooth dual there are no orbits
of codimension 0 or 1, as $a$ is a Casimir function, and in each hyperplane
the orbits are of codimension at least 1.\footnote{There exist Virasoro 
orbits 
of codimension 1 if in the dual space $vir^*$ besides smooth elements 
we also admit singular ones, cf. \cite{Wit}. In this paper we 
consider the classification of the smooth dual elements only.}

For $a_0\not=0$ the orbits  are represented
by the Hill's  operators, whose holonomy matrices were classified above,
while for $a_0=0$ they are quadratic differentials. Our choices of
representatives for the orbits of codimension 1 for a fixed $a_0$
(i.e. of total codimension 2 in $vir^*$) will be as follows.

a) For a generic point representing the case $i)$ above we take
 Hill's operator $-\partial_{x}^2+1/2\,\,(u_0=(dx)^2/2, a_0=1) $.
It corresponds to the differential equation $(\partial_{x}^2-1/2)y=0$
and has the holonomy matrix 
${\rm diag}\left(\exp(\pi\sqrt 2), \exp(-\pi\sqrt 2)\right)$,
the exponent  of
the type $i)$. This freezing point corresponds to the CH equation.

b) The  matrix of type $ii)$ can be encountered in
a generic 1-parameter family of matrices in $sl_2(\mathbb R)$. 
Its exponent, a Jordan $2\times 2$-block
with the eigenvalue 1,  can be represented as the holonomy
matrix by the Hill's operator $-\partial_{x}^2~~(u_0=0, a_0=1)$.
The latter point in $vir^*$ corresponds (after an appropriate reduction)
to the  HS equation.

c) The hyperplane $a_0=0$ in the Virasoro dual is the dual space $vect^*$
to the non-extended Lie algebra of  vector fields.
The orbits of codimension 1 in the  space $vect^*$ are represented, e.g.,  by
the quadratic differential ${\frac 12}(dx)^2$ (i.e., by
the point $(a_0=0, u_0=(dx)^2/2)$) as we discussed in Section $A$.
Freezing the Poisson structure at the latter point leads to the 
KdV equation.

The above three cases are described in Corollary \ref{classpairs}.
\end{rem}

\medskip

\begin{rem}\label{kdvint}
 Recall that the holonomy matrix $M$ of
 Hill's operator $\partial^2_x+u(x)$ changes to a conjugate one under
the Virasoro action. This implies that
$h (\partial^2_x+u(x))=\log ({\rm trace}\, M)$ is a Casimir function 
on $vir^*$.
One can use it to generate the KdV hierarchy via the Lenard scheme described in
Theorem \ref{biham}.

Recall that for the KdV equation the freezing point for the constant
Poisson structure $\{.,.\}_0$
is  $(a_0=0, u_0=(dx)^2/2)$.
Therefore, the Casimir function of the bracket
$\{.,.\}_\lambda:=\{.,.\}_{LP}-\lambda^2\{.,.\}_0$ has the form
$h_\lambda(\partial^2_x+u(x)):=\log ({\rm trace}\, M_\lambda)$,
where $M_\lambda$ is the holonomy of
the  operator $\partial^2_x+u(x)-\lambda^2$.
The expansion of the function $h_\lambda$ in $\lambda$
 produces the first integrals
of the KdV equation:
$$
 h_\lambda
\approx
{2\pi}\lambda-\sum_{n=1}^{\infty} c_n h_{2n-1} \lambda^{1-2n},
$$
where
$$
h_{1} =  \int_{S^1} u(x) \, dx,
\quad
h_{3} =  \int_{S^1} u^2(x) \, dx,  \quad h_{5}
=
 \int_{S^1} \left(
u^3(x) - \frac{1}{2} (u_{x}(x))^2 \right) \, dx, \;\;  \dots
$$
and $c_1=1/2,~~c_n=(2n-3)!!/(2^nn!)$ for $n>1$. One can see that
the Hamiltonian $h_3$ is quadratic in $u$ and coincides with the
``energy'' Hamiltonian of the KdV equation, regarded as an Euler equation.
(Note that the KdV Hamiltonians $h_j$ are differential polynomials
whose degree increases with $j$. The latter follows from
the recurrence relation 
$\{h_{2j+1},f\}_0+\{h_{2j-1},f\}_{LP}=0$ for Hamiltonians $h_j$ 
(cf. Theorem \ref{biham}) for the constant and linear Poisson brackets
on $vir^*$.) Much more details on the KdV structures can be found in \cite{gz}.

Similar computations can be done for
 CH and HS, the other two equations considered in this paper, cf. \cite{bss}.
\end{rem}

%\bigskip
%%%%%%%%%%%%%%%%%%%%%%%%%%%%%%%%%%%%%%%%%%%

\subsection*{Acknowledgements}
We are grateful for the hospitality to the  Isaac Newton Institute
(Cambridge, UK) and the Fields Institute (Toronto, Canada),
where most of this work was done.
We are indebted to P.~Pushkar and D.~Novikov for drawing  the figures.
Special thanks go to I.~Zakharevich for many useful remarks, as
well as for illuminating discussions on Theorem \ref{homEuler'}.

The work of B.K. was partially supported by an Alfred P. Sloan
Research Fellowship and by an NSERC research grant.
The work of G.M. was supported in part by NSF Grant DMS-9970857.

%\bigskip
\bigskip

%\centerline{\bf REFERENCES}
\bibliographystyle{amsplain}

\end{document}